\documentclass{article}[12pt]
\usepackage{amssymb,amsmath,amscd,amsthm}
\usepackage{graphicx,psfrag,epsfig,multirow}
\begin{document}
\title{Arnold diffusion in a pendulum lattice.}
\author{Vadim Kaloshin\thanks{partially supported by
NSF grant DMS-0701271}, \
Mark Levi\thanks{Partially supported by NSF grant DMS-0605878 }, \
Marya Saprykina }
\maketitle
\newtheorem{Thm}{Theorem}[section]
\newtheorem{Def}{Definition}
\newtheorem{thm}{Theorem}
\newtheorem{dfn}{Definition}
\newtheorem{Lm}{Lemma}
\newtheorem{lem}{Lemma}
\newtheorem{Prop}[Lm]{Proposition}
\newtheorem{Rem}[Thm]{Remark}
\newtheorem{rem}{Remark}
\newtheorem{Cor}[Lm]{Corollary}

\def\bdef{\begin{Def}}
\def\endef{\end{Def}}
\def\bthm{\begin{Thm}}
\def\ethm{\end{Thm}}
\def\bprop{\begin{Prop}}
\def\enprop{\end{Prop}}
\def\blm{\begin{Lm}}
\def\elm{\end{Lm}}
\def\bcor{\begin{Cor}}
\def\ecor{\end{Cor}}
\def\brm{\begin{Rem}}
\def\erm{\end{Rem}}
\def\bfig{\begin{picture}}
\def\efig{\end{picture}}
\def\be{\begin{eqnarray}}
\def\ee{\end{eqnarray}}
\def\beal{\begin{aligned}}
\def\enal{\end{aligned}}
\def\om{\omega}
\def\sg{\sigma}
\def\~{\tilde}
\def\Bbb{\mathbb}
\def\Lb{\Lambda}
\def\A{\Bbb A}
\def\B{\Bbb B}
\def\R{\Bbb R}
\def\C{\Bbb C}
\def\Z{\Bbb Z}
\def\T{\Bbb T}
\def\Ss{\Bbb S}
\def\cQ{\mathbb Q}
\def\cP{\mathbb P}
\def\Cal{\mathcal}
\def\cal{\mathcal}
\def\L{\mathcal L}
\def\cI{\mathcal I}
\def\dt{\delta}
\def\Si{\Sigma}
\def\bt{\beta}
\def\th{\theta}
\def\Th{\Theta}
\def\eps{\varepsilon}
\def\lb{\lambda}
\def\gm{\gamma}
\def\al{\alpha}
\def\I{\mathcal I}
\def\cT{\mathcal T}

The main model studied in this paper is a lattice of nearest neighbors
coupled pendula. For certain localized coupling we prove existence
of energy transfer and estimate its speed.  

 \section{The description of the motion.}
We consider a system of pendula with a nearest neighbors coupling:
\begin{equation}
	\ddot x_i + \sin x _i = 
	-\varepsilon \frac{\partial}{\partial x_i} \beta  (x_{i-1}, x_i, x_{i+1},\eps), \
i \in \Z,
 	\label{eq:mainsystem}
\end{equation}
where the interaction potential $\beta$ is localized and will be defined later.
This system can be written in the Hamiltonian form with 
${\bf x} =\{x_i\}_{i\in \Z},\ {\bf y} =\{y_i\}_{i\in \Z},\ x_i$ and $y_i \in \R$, with the 
the Hamiltonian
\be \label{eq:Hamiltonian}
\beal
H_\eps ({\bf x} ,{\bf y} )=\sum_{i\in \Z} \ \ \ \dfrac{y_i^2}{2}
+(-\cos x_i-1)+\eps \beta (x_{i-1}, x_i, x_{i+1},\eps)=
\\
=\sum_{i\in \Z}
\ \ \ \dfrac{y_i^2}{2}+V(x_i)+\eps \beta (x_{i-1}, x_i, x_{i+1},\eps),
\enal
\ee
where $V(x)=-\cos x-1$ is the pendulum potential.


   \begin{figure}[thb]
         \psfrag{g}[c]{\small $w$}
 	 \psfrag{w}[c]{\small $W^u$}
 	 \psfrag{-w}[c]{\small $W^s$}
    	\center{ \includegraphics[scale=0.5]{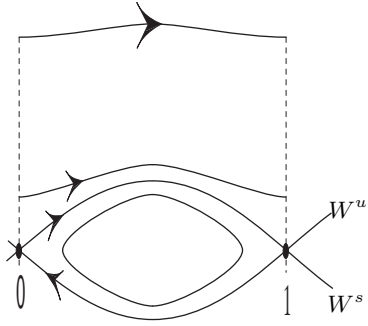}}
    	\caption{The ``running" and the near--heteroclinic motion are the building
	blocks of the dynamics.  }
    	\label{fig:pendulumbasic}
\end{figure}

\vskip 0.3 in 

The system is near--integrable for small $\varepsilon$, and most
(in the sense of measure) of the systems's phase space is 
taken up by invariant KAM tori. In particular, for most initial data
the energy of each pendulum will stay close to its initial value for all time.
Nevertheless, we will show that this is not so for some motions, where
the energy
can slowly ``seep" from one pendulum to another. We will in fact prove
that for an arbitrarily small $\varepsilon$ and for any sequence of integers
$ \sigma = (\ldots, \sigma _{-1}, \sigma_0, \sigma_1, \ldots) $ such that
$\sigma_0=0, |\sigma_j-\sigma_{j+1}|=1$ for all $j\in \Z$ there exists
a sequence of times $ (\ldots, t_{-1}, t_0,t_1,\ldots) $ (depending on 
$\varepsilon$) such that at time $ t_k $ the $\sigma_j$-th
pendulum has most   of
the system's energy. 
In particular, one can make the energy wander along
the chain of the pendula in any prescribed fashion, advancing to the right 
any number of steps, retreating to the left any number of steps, and so on.

\vskip 0.3 in 
From now on we fix the energy of the system to be $1$\footnote{Energy $0$ corresponds to 
all pendula upside down and at rest. Indeed,  the maximum of the potential energy  
$ V(x)=- \cos x - 1 $ of an individual pendulum is $0$ and is achieved at $ x= \pi $, an upside--down position.}. Below 
we shall concentrate on the case of periodic collection of 4 pendula, 
i.e. of the index $i$  (mod $4$). The proof in the general periodic case 
$ i \in {\mathbb Z} /p {\mathbb Z}$
is quite similar and necessary remarks are made along 
the proof.

    \begin{figure}[thb]
         \psfrag{g=e}[c]{\small $g    $}
 	 \psfrag{x1}[c]{\small $x_1$}
	  \psfrag{x2}[c]{\small $x_2$}
	  \psfrag{x3}[c]{\small $x_3$}
	  \psfrag{x4}[c]{\small $x_4$}
	  \psfrag{g}[c]{\small $w$}
	  \psfrag{g}[c]{\small $w$}
	  \psfrag{g}[c]{\small $w$}
    	\center{ \includegraphics[scale=0.7]{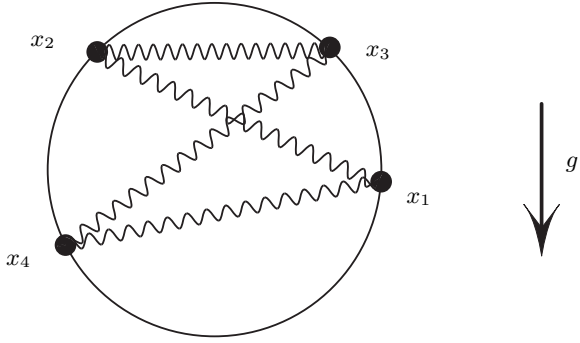}}
    	\caption{A mechanical interpretation of (\ref{eq:sinegordon}) with $ \beta =
	  \sin (x_{j-1}-2 x_j+x_{j+1} ) $. }
    	\label{fig:pointsoncircle}
\end{figure}


We note as a side remark that the space discretization of the sin-Gordon equation
$ u_{tt}-u_{ss}= \sin u $ 
results in a system of pendula with elastic coupling \cite{BWE, KP, WES}:
\begin{equation}
  	 \beta (x_{j-1}, x_j, x_{j+1}) = a(x_{j-1}+x_{j+1}-2x_j);
	 \label{eq:sinegordon}
\end{equation}
 this corresponds to an elastic torsional coupling between the neighbors. 
 In particular, as
 the angle $ x_{j+1}-x_j \rightarrow \infty $ we have $ \beta \rightarrow \infty $.
 By contrast, the coupling we consider in this paper can be interpreted as coming
 from a spring connecting points on a circle with angular coordinates $ x_j $, 
 as shown in Figure~\ref{fig:pointsoncircle}.

 The coupling in our system is, however,  localized, as described next. 
The class of coupling functions $\beta$ for which our results   hold is defined as follows: 
Let $\eta:\R_+\to \R$ be
a $C^\infty$ bump function: $\eta(x)>0$ for $|x|<1$ and $\eta(x)=0$ for $|x|\ge 1$. 
Exact form of $\eta$ is not important, and in particular, no monotonicity properties are assumed; in particular, $\eta $ is allowed to have many local maxima and minima, as long as the above conditions
hold. We now define
\begin{equation} 
\beta(x,\eps)=\eps^r \sum_{n \in \Z^3} \eta \left(
\dfrac{|x - 2 \pi n|}{\eps}\right), \qquad x=(x_1,x_2,x_3)\in \R^3, 
\label{eq:beta}
\end{equation} 
and fix $ r \geq 3 $ from now on. 
This is a $C^\infty$-smooth $2\pi$-periodic function in each $x_j, \, j=1,2,3,$
or, equivalently, a function on $2\pi \R^3 / \Z^3$. Note that the $C^r$-norm
of $\eps \beta(\cdot , \eps)$ tends to zero as $\eps \to 0$, while the norms of
order $r+2$ and higher are unbounded for $ \varepsilon \rightarrow 0 $.  
We fix any finite $r\ge 3$ from now on.

We will sometimes refer to connected components of support of $\beta$ as {\it lenses}:  in fact, they 
act by defocusing geodesics in the Jacobi metric, as explained later, as in 
\cite{KL1}. 
 
\vskip 0.3 in  

\begin{figure}[thb]
         \psfrag{g}[c]{\small $w$}
	\center{  \includegraphics[scale=0.7]{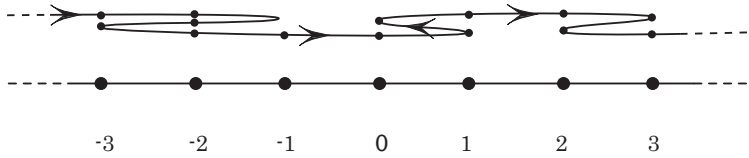}}
	\caption{Any path in the graph $ {\mathbb Z}  $  can be shadowed by a solution of   (\ref{eq:mainsystem}). }
	\label{fig:graph}
\end{figure}

\noindent According to the  main theorem, stated next,   the energy 
\begin{equation} \label{eq:energy}
	E_j:=\dfrac{\dot x _j^2}{2} + V(x_j) 
\end{equation}
at the $j$th site can  pass  from one site to another  in an arbitrarily prescribed  sequence of steps, as illustrated in Figure~\ref{fig:graph}. Here is a more precise statement. 
\begin{thm} \label{thm:main}
Let us fix the total energy\footnote{In fact, any value in excess 
of the  potential energy of an upside--down equilibrium works.} $ E=1$ in 
the system  (\ref{eq:mainsystem}) with $\beta$ satisfying (\ref{eq:beta}). 
There exists $ \varepsilon _0 > 0 $ such that for any 
$ 0 < \varepsilon < \varepsilon_0 $ and for any path 
$ \ldots \sigma _{-1} \sigma _0 \sigma _1 \ldots $ in the graph $ {\mathbb Z}  $ 
there exists a solution of   (\ref{eq:mainsystem}) and a sequence of times 
$ \ldots t_{-1}t_0t_1\ldots $ such that  the energies   (\ref{eq:energy})  of 
individual pendula satisfy
$$
	| E_{\sigma_j}(t_j) - 1| < C  \sqrt{ \varepsilon }, \ \ \hbox{and}\ \ \ 
	| E_\sigma(t_j) | < C  \sqrt { \varepsilon } \ \ \hbox{for}\ \ \  
	\sigma \not=  \sigma_j, 
$$
where $C$ is independent of $\varepsilon$. The times $ t_j $ can be chosen so that 
\begin{equation} 
0< t_{j+1}-t_j\leq C \varepsilon^{-4r-8}. 
 \label{eq:steptime} 
\end{equation} 
\end{thm}

This theorem shows that, although the system   (\ref{eq:mainsystem})  is 
near--integrable, so that for most (in the sense of Liouville measure) solutions 
the action stays close to its initial value for all time,  there exist nevertheless solutions 
for which the action changes by $ O(1) $ no matter how small $\varepsilon$ is. 
In other words, the system exhibits Arnold diffusion. According to (\ref{eq:steptime})  
the  rate of this diffusion   is polynomial. The bound in  (\ref{eq:steptime})  is 
not sharp, but it can be improved by a more careful tracing of the estimates 
in our example. In the general case, polynomial upper
bounds for speed of diffusion for finitely-differentiable systems 
have been obtained in \cite{Bu}.

The first example of Arnold diffusion was outlined in the well known  paper of Arnold \cite{A}.
Bessi \cite{Bs} (see also \cite{BB}) proved diffusion in Arnold's example by a variational method, by considering the gradient flow of the Lagrangian action functional. John Mather \cite{Ma} used a somewhat similar approach 
to construct accelerating orbits for time periodic mechanical systems on a $2$-torus
(see also \cite{BT, GT, DLS, Ka}). References to the recent progress on Arnold diffusion goes 
beyond the scope of this paper and can be found in \cite{KL2}. 
In the present paper we use a slightly different version  
of this approach, based on using the Maupertuis' principle. We construct the  ``diffusing" solutions as  geodesics in a Jacobi metric, so that all these solutions have a fixed prescribed energy.  These geodesics are constructed by concatenating geodesic segments which follow a prescribed itinerary. 
The construction is fairly similar to 
\cite{KL1,KL2}. 

Anderson localization is an important example of energy 
(non)transfer (see \cite{LTW} for a survey), still not very well understood. The role of Arnold 
diffusion for destruction of Anderson localization is discussed in \cite{Ba}. 
Probably the most popularized   lattice model is the one introduced by Fermi-Pasta-Ulam 
in their seminal paper \cite{FPU}. Although most small amplitude solutions in the FPU model
do not exhibit energy transfer (see e.g. \cite{HK}), proving the existence 
of solutions with energy transfer is an interesting open problem. 
Other physically significant lattice models are discussed in \cite{FSW}. 

Understanding of the transfer of energy for Hamiltonian PDEs is one 
of emerging directions of research  (see \cite{Bo}; a recent progress for 
the cubic defocusing nonlinear Schr\"odinger equation has been made in \cite{I}).

\begin{figure}[thb]
	   \psfrag{s1}[c]{\ }
	   \psfrag{s2}[c]{\ }
	   \psfrag{s3}[c]{\ }
    	\center{ \includegraphics[scale=0.7]{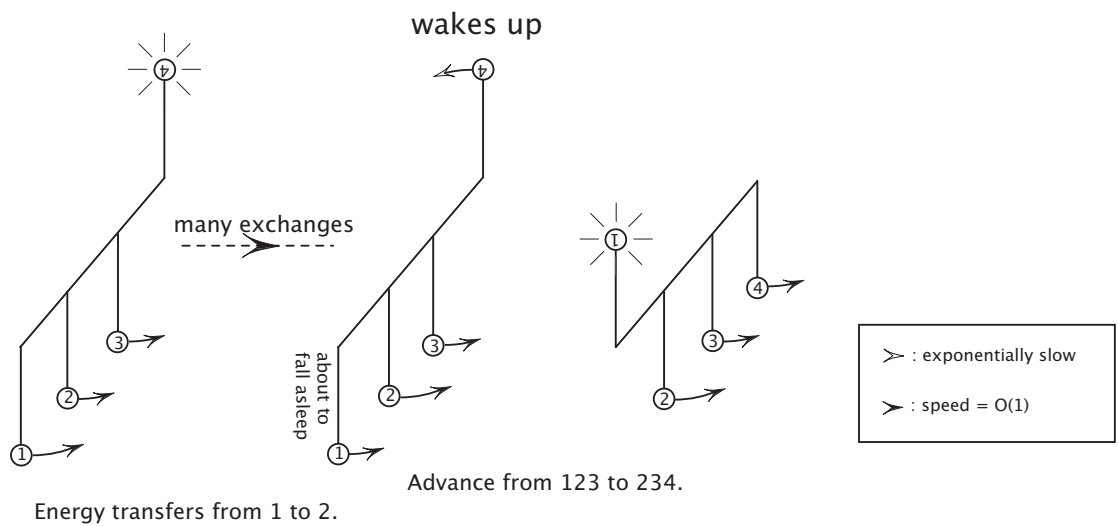}}
    	\caption{One full step in the propagation of the ``kink".    }
    	\label{fig:pendula}
\end{figure}

\subsection{A heuristic description of propagation}

In this section we give a purely heuristic picture of the physical 
motions showing Arnold diffusion. As mentioned earlier, we consider the periodic case
$ x_{i+4}=x_i $ as a representative example. 

\vskip 0.1 in 

\begin{figure}[thb]
	  \psfrag{s0}[c]{\small $\Sigma_{12,3}^1$}
	   \psfrag{s1}[c]{\small $\Sigma_{12,3}^{N_1}$}
	  \psfrag{s2}[c]{\small $\Sigma_{23,4}^1$ }
	  \psfrag{s3}[c]{\small $\Sigma_{23,4}^{N_2}$  }
   	\center{ \includegraphics[scale=0.3]{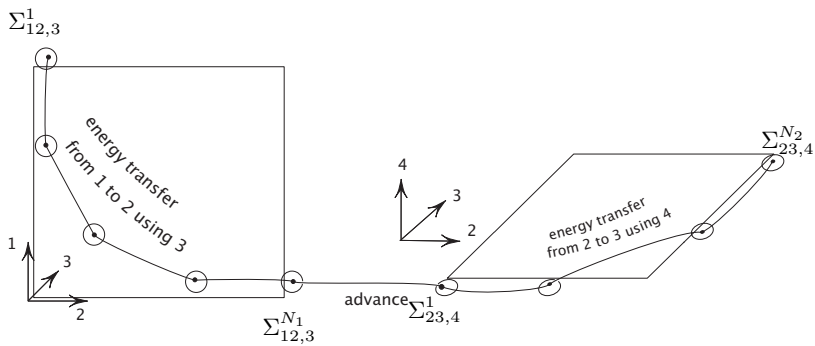}}
	\label{fig:configspace}
    	\caption{Energy transfer and sections in the configuration space $ {\mathbb R} ^4 $.}
	\label{fig:sections}
\end{figure}

\paragraph {Stage 1:  transfer of energy.} In this stage only three pendula:
{\bf 1}, {\bf 2}  and {\bf 3} governed by 
(\ref{eq:mainsystem})  with $ k = 1, 2 $ and $3$ are ``active", while {\bf 4} 
``sleeps" upside-down (see Figure~\ref{fig:pendula}, left). The stage consists 
of many substages  illustrated by Figures~\ref{fig:sections}, left.  In each of 
these substages a small amount of energy is transferred from {\bf 1} to {\bf 2}. 
This transfer is somewhat similar to the one described in \cite{KL1} for a metric 
on the $3$--torus. 
Finally, in the last substage,  {\bf 1}  is left with just enough energy to climb 
upside down and to fall asleep there, while {\bf 2} rotates with speed $ O(1) $,
as shown in the middle of  Figure~\ref{fig:pendula}. The same motion, viewed in the 
configuration space $ {\mathbb R}  ^4  $, is shown in Figure~\ref{fig:sections}. 

The motion just described is similar to that in a slightly simpler example described in \cite{KL1, KL2}.



   \paragraph {Stage 2:  advance.} This stage is sketched in
Figure~\ref{fig:pendula}, right and Figure~\ref{fig:sections}, middle.  
At $ t=t_1 $ three neighbors, say {\bf 1}, {\bf 2},
and {\bf 3},   are in the bottom position. The middle pendulum {\bf 2}
is running: $ \dot x _2=O(1) $, while its two neighbors {\bf 1}   and {\bf 3}
have near--heteroclinic  speeds close to the heteroclinic speeds
$\sqrt {-2 V(x_i)},\ i=1,3$ resp., at $ t=t_1 $. The remaining pendulum {\bf 4}
is up (see Figure \ref{fig:pendula} middle). As the time goes on, while {\bf 2}
is spinning
with speed $ O(1) $,   {\bf 1} rises to the top equilibrium, where it will sleep until
further notice, while the sleeper {\bf 4}  ``wakes up", i.e. falls from its perch,
turning by $\pi$ by the exact moment when 
 {\bf 2} finishes a large integer number  of full  spins. By that moment,
 $ x_3 $ makes a ``gentle" turn by $ 2 \pi $, returning to the bottom
position. In short, the accomplishment of this stage is 
the falling asleep of\ {\bf 1} and the awakening of \ {\bf 4}. This is illustrated 
in Figure~\ref{fig:pendula}, right.  We will call this stage the ``advance" 
because of its similarity with the advancing caterpillar: a rear foot {\bf 1} 
is placed on the ground, while the front foot {\bf 4} is lifted, ready to move.

    The ending moment of the second stage is the beginning moment of
   the first stage described above modulo the shift of the index by $1$.
   We have, in other words, a ``traveling wave" -- a (very) discrete analog of
   the kink in the sine--Gordon equation. However, in contrast to the standard
   traveling kink, ours can change the direction of its propagation arbitrarily,
   according to a prescribed itinerary.

\section {The  proof of Theorem \ref{thm:main}. }   \label{sec:proof}

The full complexity of the problem is already seen in 
the case of $4$ pendula, and we limit our consideration to that case. 
Now we restate the theorem in geometrical terms. The following is motivated by 
the heuristic outline of the energy transfer between the pendula: as 
mentioned before (see Figures~\ref{fig:pendula} and~\ref{fig:sections}), we want the energy  
to pass from one pendulum (e.g. {\bf 1}) to another (e.g. {\bf 2})
in small increments over many steps, after which the pendula should 
change roles: the ``giver" {\bf 1} ``falls asleep", while 
a ``sleeper" (e.g. {\bf 4}) ``wakes up" thus enabling the ``taker" {\bf 2}
to become the next ``giver".  This is reflected in the following 
geometrical construction. First, we construct an itinerary 
(section \ref{itinerary}) for the desired orbit. Then we reformulate 
the problem of existence such an orbit as a variational problem.  
Finally, we prove the existence of a solution to this problem in 
section \ref{interior-minimum}. In the latter section we use 
two lemmas \ref{lem:existence} and \ref{lem:velocity} stated (and proved) in 
sections \ref{sec:existence} and \ref{sec:hyperbolic} respectively.

\subsection{Constructing an itinerary}\label{itinerary}
  \begin{figure}[thb]
         \psfrag{B}[c]{\small $\Sigma_{12,3}^{N_1}$}
 	\psfrag{B1}[c]{\small $\Sigma_{23,4}^1$}
	\center{  \includegraphics[scale=0.6]{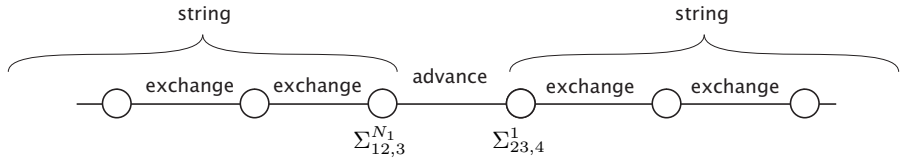}}
	\caption{An itinerary. }
	\label{fig:groups}
\end{figure} 

To prove Theorem \ref{thm:main}, given a path in the graph 
(see Figure \ref{fig:graph}) we must produce an orbit that shadows this path. 
We describe the construction in the case of {\it monotone 
energy transfer} (i.e., to the right neighbor):  $\sigma_{j+1}=\sigma_j + 1=j+1$ for 
all $j\in \Z$. Extension to the case of general path 
$(\dots \sg_0 \sg_1 \sg_2 \dots)$ poses no new difficulties. 
We this consider an infinite  sequence of codimension one sections in 
${\mathbb R}^4$,  grouped into finite strings, Figures~\ref{fig:sections} 
and~\ref{fig:groups}: 
\begin{equation} 
	\ldots  \ 
	\underbrace{(\Sigma_{12,3}^1  \Sigma_{12,3}^2  
	\ldots  \Sigma_{12,3}^{N_1})}_{1\to 2} \ 
	\underbrace{(\Sigma_{23,4}^1, \Sigma_{23,4}^2, 
	\ldots, \Sigma_{23,4}^{N_2})}_{2\to 3} \ldots  , 
	\label{eq:sequencofballs}
\end{equation} 
where the sections and their spacings are defined according to the following rules:  
\begin{enumerate} 
\item Section  $ \Sigma_{12,3}^1 $, for example, is seen in 
Figure~\ref{fig:pendula}, left. The subscripts $ 12$ indicate  
that  $1$ and $2$ exchange energy, and $3$ is the ``facilitator":     
\begin{equation} 
\Sigma_{12,3}^1:=\{  x_3=0, \ \   x_1 ^2 + x_2 ^2 \leq \varepsilon,
 \     |x_4- \pi | \leq   \sqrt \varepsilon  \}.
\label{eq:section1}
\end{equation} 

\item 
All   sections within each string in   (\ref{eq:sequencofballs})  are integer translates 
of each other. For example, in the first string $ 1\to 2 $: 
$$
 	\Sigma_{12,3}^{k+1} = \Sigma_{12,3}^{k} + 
	2 \pi (m_{12,3}^k,n_{12,3}^k, 1,0)=:
	\Sigma_{12,3}^{k}  +\vec n_{12,3}^{k} \qquad k=1, \dots, N_1
$$
To define the second string, we  replace $i \mapsto i+1$ mod $4$
in  (\ref{eq:section1}), setting
\begin{equation} 
\Sigma_{23,4}:=\{  x_4=0, \ \   x_2 ^2 + x_3 ^2 \leq \varepsilon,
 \     |x_1 - \pi | \leq   \sqrt \varepsilon  \}.
\label{eq:section2}
\end{equation} 
We then define
$\Sigma_{23,4}^1$ as a translate of $\Sigma_{23,4}$ by  integer multiples  
of $2\pi$ in each coordinate.
Each section  in the second  string is an  integer translate 
of the previous section   by  
\[
\vec n_{23,4}^{k}:=2 \pi (0,m_{23,4}^k,n_{23,4}^k, 1)
\] 

\item The neighboring strings are related via 
$$
Center( \Sigma^{1 }_{ 23,4})=Center(\Sigma^{N_1}_{12,3}) +
 2 \pi ( \frac{1}{2},   m_{23,4}^0,  1, \frac{1}{2})=
 Center(\Sigma^{N_1}_{12,3}) + \vec  n_{23,4}^0.
$$
 
\item All translates are far apart: each 
\begin{equation} 
	 | \vec n^k_{\cdots} | \geq \varepsilon ^{-2r-4} 
	\label{eq:long}
\end{equation}   

\item The turns are gradual   in the sense that the unit vectors 
	$ e^k_{\cdots}= \vec n^k_{\cdots}/ | \vec n^k_{\cdots} | $ satisfy 
\begin{equation} 
	| e^{k+1}_{\cdots}-e^k_{\cdots} | \leq \varepsilon ^{2r+4} , \ \ \
	 |  e^1_{\cdots+1}-e^{N_{j}}_{\cdots} | \leq \varepsilon ^{2r+4}. 
	\label{eq:slowturn} 
\end{equation}  
\end{enumerate} 
As mentioned above, treating a   general itinerary poses no difficulties. 


 \subsection{A variational problem and its solution} \label{interior-minimum}
 We note that the energy one solution of   (\ref{eq:mainsystem})  are
the geodesics in the Jacobi metric\footnote{Up to the factor $  \sqrt{ 2 } $, the square root is the 
speed of the energy one solution in the configuration space.} 
\begin{equation} 
	d\rho({\bf x} )= \sqrt { 1-\sum_{i=1}^4 \large(V(x_i)- 
	\varepsilon \beta (x_{i-1},x_i,x_{i+1} , \varepsilon )\large)} \ \ ds, 
	\ \ \ 
	\label{eq:jacobimetric}
\end{equation} 
where ${\bf x} =(x_1,x_2,x_3,x_4), \ V(x)=-\cos x-1\le 0$, 
$ ds$ is the Euclidean metric and indices of $x_i$'s are taken mod $4$. 

With the sections  having been defined  in items 1--5 above, 
we now list {\it  the main steps of the    proof  of Theorem \ref{thm:main},}  
and fill in the details in the following sections.  
  \begin{enumerate}
\item {\it Defining geodesic segments.} Let $ \Sigma_0 , \Sigma _1 $ be
two consecutive sections in the chain of sections   (\ref{eq:sequencofballs})    
and let  $ p_i\in\Sigma_1 $, $ i=0, 1 $. Centers of these sections
differ by $2\pi \vec n$ with $\vec n$ being either
$(m,n,1,0)\in \Z^4$, or $(\frac 12, s, 1, \frac 12)$ with $s\in \Z$ and 
satisfying (\ref{eq:long}). According to Lemma~\ref{lem:existence}
from section \ref{sec:existence}, there exists a connecting
geodesic $ \gamma (p_0,p_1) $ of (\ref{eq:jacobimetric}) 
  which depends
smoothly on its ends $ p_0 ,\  p_1 $. At this stage the integer
parameters either $m$ and  $n$ or $s$ are still free.

\item {\it Constructing a long shadowing geodesic.} Consider 
 a finite  segment of $N+2$ sections from the sequence (\ref{eq:sequencofballs}).  To simplify notation, we denote these sections by  $ \Sigma _i $, $ 0\leq i \leq N+1 $ ($N$ here is arbitrarily
large). We also choose arbitrary points $ p_i\in \Sigma_i $. Later we will treat $ p_0, p_{N+1}$ 
as fixed and $ p_1, \ldots , p_N $ as variable. 

According to the preceding item there exists a broken geodesic 
$ \gamma (p_0, \ldots, p_N ) $, a concatenation of energy one orbits $ \gamma (p_i, p_{i+1})$ of   (\ref{eq:mainsystem}).

The   length (in the Jacobi metric) of this broken geodesic
\begin{equation}
 	 L(p_1, \ldots, p_N) = L(p_0,p_1)+\ldots + L(p_{N-1},p_ N).
	 \label{eq:sumsegments}
\end{equation}
 is a function of break points $ p_j$; we omit $ p_0 $ and $ p_{N+1} $ from 
 the left--hand side since they will be  considered as fixed.
 


We will  show that $ L $ has {\bf a minimum on the interior of its domain}
$ \Sigma_1 \times \ldots \times  \Sigma_N $.    Such an interior minimum corresponds  to a true geodesic. We will thus establish the existence of a geodesic with a prescribed itinerary. 

\item {\it The key step: existence of an interior minimum for   (\ref{eq:sumsegments}).   }
Let us consider  two consecutive terms from the sum (\ref{eq:sumsegments}):
\begin{equation}
	S(p)= L(p_{j-1},p )+L(p,p_{j+1}), \ p \in \Sigma_i
	\label{eq:twosegments}
\end{equation}
where $ p_{j\pm 1}\in \Sigma_{j\pm 1} $   are fixed and $ p\in \Sigma_i $
is  variable. Without loss of generality we take $\Sigma_j=\Sigma_{12,3}^1$,
given by (\ref{eq:section1})\footnote{For the future reference, we note that 
the triple may or may not be entirely 
in one string in the sequence    (\ref{eq:sequencofballs})  of sections. }.
To prove the existence of an interior minimum of   (\ref{eq:sumsegments}) it suffices 
to show that the minimum of $ S(p) $ is achieved in the interior of $ \Sigma_j $. 
To that end we first alter $\beta$   in the Jacobi metric  (\ref{eq:jacobimetric}) by setting $ \beta = 0 $ {\it  only }  in the cylinder $ x_1 ^2 + x_2 ^2 + x_3 ^2 < \varepsilon ^2 $ which passes through the center of  $ \Sigma_j $ (we do not alter $\beta$ anywhere else).  Before restoring $\beta$ to its original form, we study the associated length 
\be \label{eq:action}
	S^0(p)= L^0(p_0,p )+L^0(p,p_2), \ p \in \Sigma_j\equiv\Sigma.  
\ee
 Once the properties of   $ S^0 $ are established 
(see   (\ref{eq:flat})  and   (\ref{eq:convex}) below), we will show that restoring $\beta$ to its original form creates a minimum for $ S $ in the interior of $ \Sigma $. Without the loss of generality, we take $ \Sigma = \Sigma_{12,3}^1 $ as in   (\ref{eq:section1}), so that 
$ S^0=S^0(x_1,x_2,x_4) $. 
  
We   will first show 
  that $ S_0 = S^0(x_1,x_2, x_4) $ is ``nearly constant" in the first two variables and has
a minimum near the ``equator" $ x_4= \pi $: 
 
 \begin{equation}  
| S^0(x_1,x_2,x_4) - S^0(0,0,x_4)| \le 2 c  \eps^{r+2.5},
	\label{eq:flat}
\end{equation}   
\begin{equation} 
	S^0(x_1,x_2, \pi\pm  \sqrt{ \varepsilon } )> S^0(x_1,x_2,\pi)+\frac{\varepsilon }{2} 
	\label{eq:convex}
\end{equation}      
	for any $(x_1,x_2,\pi)\in \Sigma  $.  
\end{enumerate}

\paragraph {Proof of   (\ref{eq:flat}) and   (\ref{eq:convex}).} 
By Lemma 1
from \cite{KL1} we have, with $p=(x_1,x_2,x_3,x_4)$:
$$
\frac{\partial L^0(p_-,p)}{\partial x_i} = \dot x^-_i,
 \ \ \frac{\partial L^0(p ,p_+)}{\partial x_i} = -\dot x^+_i, 
$$
where $ {\bf x} ^- (t) =(x^-_1, x^-_2, x^-_3, x^-_4)$ is the energy one solution   with the 
modified $\beta$, 
connecting $ p_{j-1} $ to $ p $, and where the differentiation is taken
at the moment the solution passes through $p$. This solution exists
by Lemma \ref{lem:existence} below.    We use a similar notation  $ x^+$
for the energy one solution connecting $ p $ with $ p_{j+1} $.  
We thus conclude that
\begin{equation}
	\frac{\partial S^0}{\partial x_i}(p) = \dot x^-_i-\dot x^+_i , \ \ i=1, 2, 4; 
		\label{eq:gradS}
\end{equation}
this identity\footnote{we do not differentiate by $ x_3 $ since we only need to define 
$ S $ and $ S^0 $  on $ \Sigma \subset \{x_3=0\}$.} will allow us to analyze $S^0$.  Now, due to the fact
that the perturbation $\beta$ near $p$ is removed the pendula are decoupled
and the velocity is
explicitly given in terms of energy distribution (\ref{eq:energy})
$$
|\dot x_i| =   \sqrt{ 2(E_i-V(x_i) )},
$$
where $ E_i $ is the energy of the $i$th pendulum near $p$
The estimate of  (\ref{eq:gradS})  is now reduced  to studying the difference of velocities.  We have to consider  two    cases: in one, $ \Sigma_{j-1}, \ \Sigma, \ \Sigma_{j+1}$     belong to the same string in   (\ref{eq:sequencofballs}) (the case of ``energy transfer''), and in the other,    they do not (``the advance").

\paragraph {Case 1 -- energy transfer.  }
In this case   all sections lie in the same string in   (\ref{eq:sequencofballs}) --  say, in  
$ 1\to 2 $.   The displacements  $2\pi \vec n_-$ and $2\pi \vec n_+$
are then of the form $2\pi(m_\pm,n_\pm,1,0)$ with integer
$m_\pm$ and $n_\pm$. Assuming the integers to be positive (we can always assume them to be of  the same sign),   
   we have
\begin{equation} 
	\dot x^-_i>0,  \ \ \dot x ^+_i>0, \ \ \ i=1,2,3, 
	  \label{eq:signderiv}
\end{equation} 
at the moment when $ \Sigma $ is crossed. 


The fact that the signs are the same for $ x^+ $
and $ x^- $ is of key importance because it provides a
near--cancellation in   (\ref{eq:gradS}) for $ i=1,2,3 $.
Thus, we have
\begin{equation}
	\frac{\partial S^0}{\partial x_i} =    \sqrt{ 2(E^-_i - V(x_i)) }-
\sqrt{ 2(E^+_i - V(x_i)) }, \ \ i=1,2.
	\label{eq:S123}
\end{equation}
Note that if $ 0\leq A \leq B $ then
$ \sqrt{ B } -  \sqrt{ A } \leq  \sqrt{ B-A } $;
this, used in   (\ref{eq:S123}), gives
\begin{equation}
	\biggl|\frac{\partial S^0}{\partial x_i} \biggl|\leq
	\sqrt{ 2| E^+_i-E^-_i |  }, \ \ i=1,2.
	\label{eq:S123a}
\end{equation}

 Now according to    Lemma \ref{lem:velocity}
 and the assumptions   (\ref{eq:long})  and   (\ref{eq:slowturn}) we have  
\begin{equation}
	| E^+_i-E^-_i | \leq 	c \eps^{2k+4}, \ i=1,2,3,4
	  \label{eq:HL}
\end{equation}
and by   (\ref{eq:S123a})
\begin{equation}
	| \partial_{x_i} S^0| \leq c  \eps^{k+2}, \ i=1,2.
	\label{eq:Sestimate}
\end{equation}
Combining this with (\ref{eq:S123a}), we obtain
  (\ref{eq:flat}),   
showing that $S^0$ is ``flat" in the first two variables.

To estimate $\partial S^0/\partial x_4$ note that
$(x_4^-,\dot x_4^-)$ and $(x_4^+,\dot x_4^+)$
stay at the $\sqrt \eps$-neighborhood of
the saddle $(x_4,\dot x_4)=(\pi,0)$.
Since the distance between sections is large  (\ref{eq:long}) and the phase velocity is bounded 
(by the choice of fixed energy),    
the duration of each stage is at least $\eps^{-2k-4}$.
This implies that   $(x^-_4,\dot x^-_4)$ is at worst $O(\exp(-\eps^{-1}))$-close
to the unstable manifold 
$$
	 y = U(x)=(x- \pi )+ O((x- \pi )^2 ),
$$
 while 
$(x^+_4,\dot x^+_4)$ is at worst $O(\exp(-\eps^{-1}))$-close
to the  stable manifold $ y =- U(x) $ of the same saddle -- all this at the moment  when the solution crosses $\Sigma$. That is, 
\begin{equation} 
	\dot x_4^-=U(x)+O(\exp(\eps^{-1}))  
	\label{eq:proximity+}
\end{equation}   
 and 
\begin{equation} 
	\dot x_4^+=-U(x)+O(\exp(\eps^{-1}))  
	\label{eq:proximity-}
\end{equation}   
so that
\begin{equation}
\frac{\partial S^0}{\partial x_4} =2U(x)=
2(x_4 -\pi )+O((x_4-\pi)^2)+O(\exp(\eps^{-1})).
	\label{eq:S123b}
\end{equation}
Integration by $ x_4 $ gives    (\ref{eq:convex}).  

 

\vskip 0.3 in
\paragraph {Case 2 --- advance.} In this case not all sections lie in the same string in 
  (\ref{eq:sequencofballs}); without the loss of generality, assume that $ \Sigma_{j-1}, 
  \Sigma $ are the last two in the string $ 1 \to 2 $, while $ \Sigma_{j+1} $ is the first in the following string $ 2 \to 3 $. In this case the displacement vectors are of the form 
  $\vec n_-= (m,n,1,0)$,  $ \vec n_+=( \frac{1}{2}, s,1, \frac{1}{2} )$. In this case we still have    (\ref{eq:signderiv}), and 
 following (\ref{eq:S123}) and  (\ref{eq:S123a}) we obtain     (\ref{eq:flat}). Since the sign of 
  $ \dot x _4 $ is unknown, we treat it separately, observing, as before, that
    (\ref{eq:proximity+})  and   (\ref{eq:proximity-})  hold. This implies   (\ref{eq:S123b})  and thus       (\ref{eq:convex}).

This completes the proof of    (\ref{eq:flat})  and   (\ref{eq:convex})  in both cases. 



\paragraph {Proof of the interior minimium for $ S $.} Using the properties  of $ S^0 $ and the positivity of $\beta$ we now  show that
$S$ has a minimum inside $ \Sigma $. We do so for Case 1;
the remaining case is treated almost verbatim.  
%


 The boundary
$\partial \Sigma=\partial_v \Sigma\cup \partial_h \Sigma$ consists of
the ``vertical" and the ``horizontal" parts (after possible reindexing of the
coordinates):
\be
\beal
	 \partial_v \Sigma=\{ x_1 ^2 + x_2 ^2= \varepsilon,
	 \ | x_4 -\pi | \leq  \sqrt \varepsilon  \},\\
	 \ \ \ \partial_h \Sigma = \{ x_1 ^2 + x_2 ^2 \leq \varepsilon,
	 \  | x_4 -\pi | =  \sqrt \varepsilon  \}.
\enal
\label{eq:boundary}
\ee
A key observation we will use shortly is this:
\begin{equation}
	S(p)=S^0(p) \ \ \hbox{for all } \ \ p\in \partial_v\Sigma. 
	\label{eq:boundaryeq}
\end{equation}
Proof:    We wish to show that the energy one solution $ \gamma (p_{j-1},p) $
with  $p$ lying on $ \{x_1 ^2 + x_2 ^2 = \varepsilon , \ x_3=0\} $ does not intersect
the lens $ x_1 ^2 + x_2 ^2 + x_3 ^2 \leq \varepsilon^2 $. To that end assume
the contrary: the solution   travels from   one set to the other, taking some
time $ t = t^*>0 $. Since the distance between the sets is   
$\geq  \frac{1}{2}  \sqrt{ \varepsilon } $,
while the speed is  $\leq 2$, the time of travel is    $t^*> \frac{1}{4}  \sqrt{ \varepsilon } $. 
  But  $ \dot x_3 \geq 1 $ for
as long as $ | x_3 | \leq  \sqrt{ \varepsilon } $. Thus during   time $t^* $,   $ x_3$ changes by the amount  $ \Delta x_3 > \frac{1}{4}  \sqrt{ \varepsilon } $, which means that the solution
lies outside the lens $ x_1 ^2 + x_2 ^2 + x_3 ^2 \leq \varepsilon $,
in contradiction with the definition of $ t^* $. This proves
(\ref{eq:boundaryeq}). 

  We now show that the restriction of $S$  to each horizontal disk in $ \Sigma$: 
 $$
 	D_{h}:=  \{x _1 ^2 + x_2 ^2  \leq  \varepsilon,\  x_4= \pi +h \},\    
	   \ | h | \leq \sqrt \varepsilon 
$$
 has a minimum in the interior of $D_h $.  
To that end, we first note a crucial fact that $\beta$  
decreases $ S $ (as compared to $ S^0$) near the center
$ C_h=(0,0, \pi+h ) $ of each  $D_h $.
Note that by the definition 
(\ref{eq:beta}) the infimum of $\bt(\cdot,\eps)$ taken 
over the set $ x_1 ^2 + x_2 ^2 +x_3 ^2  \leq \varepsilon ^2 /2 $ 
is bounded from below by $b \eps^r$ for some $b>0$ 
independent of $\varepsilon$. Therefore, comparing 
the geodesic  length in the original Jacobi metric with 
the truncated one,  we obtain 
\begin{equation}
	S(C_h) \leq S^0(C_h)- \varepsilon \inf \bt(\cdot,\eps)
     \leq S^0(C)- b\varepsilon^{r + 1}.
	\label{eq:Scenter}
\end{equation}
See  \cite{KL1} (proof of Lemma 4) for more details on this argument.  

On the other hand, by   (\ref{eq:flat})  we have, for any 
$p\in \partial D_h $, $ | h | \leq  \sqrt{ \varepsilon } $: 
$$
 S^0(C_h)\leq S^0(p) + c \varepsilon ^{r+2.5} .
$$
  Combining this with   (\ref{eq:Scenter}) and   (\ref{eq:boundaryeq}) we obtain
  $$
  	S(C_h)\leq S(p)+ c  \varepsilon ^{r+2.5}
	- b\varepsilon^{r + 1}< S(p),   \ \ \forall p \in \partial D_h. 
  $$
  We showed that the minimum of $S$ cannot be achieved on $ \partial_v \Sigma $,
  and it remains to show that it cannot be achieved on $  \partial_h \Sigma $ either.
Estimate   (\ref{eq:convex}) shows that $ S^0 $ has a pronounced minimum near the equator
$ x_4= \pi $.   
   By the same estimate as we used for   (\ref{eq:Scenter}), we have   a two-sided result:
  $ | S^0(x)-S(x) | \leq \varepsilon ^{r +2} $, which together with   (\ref{eq:convex})  gives
$$
  	S(x_1, x_2,  \pi  \pm\sqrt{\varepsilon } ) > S(x_1, x_2,  \pi    ).
$$
This proves that the minimum is achieved {\it  inside } $\Sigma$.

To complete the proof of the main theorem  it remains to observe that the existence of the    the internal minimum   for $ S $ implies the existence of the internal minimum  for 
$ L(p_1, \ldots , p_N) $, as well as the existence of an internal minimum for an   infinitely long sequence   (\ref{eq:sequencofballs}). The details can be found in   \cite{L}.

 \section{The pendulum Lemma.}
 In this section we state and prove an auxiliary lemma which is used in the proofs of the main  two lemmas in the following two sections. 

   \begin{lem} \label{lem:pendulumlemma}
{\it  For any  $ |\alpha | \leq 1 $, $ | \beta - \pi | <1 $, and $T>0$
there exists a unique solution $ x(t; T, \alpha , \beta ) $ of
$ \ddot x + \sin x = 0 $, satisfying $ x(0) = \alpha $,
 $ x( T) = \beta $ with the additional property
 \begin{equation}
 	 \al \leq x(t)\leq \max\{\beta,\pi\} \ \ \hbox{for} \ \ 0\leq t\leq T.
	 \label{eq:between}
\end{equation}
This solution depends smoothly on $T$, $\alpha$ and $\beta$, and, moreover,
as $T$ increases from $0$ to $ \infty $, the energy $ E = \dot x ^2 / 2 + (- 1- \cos x) $
decreases monotonically  from   $ \infty $ to $ 0  $.

    \noindent  A similar statement holds if either
$$
	| \alpha + \pi | \leq 1, \ | \beta | \leq 1,
$$
$$
	| \alpha   | \leq 1, \ | \beta - 2 \pi  | \leq 1,
$$
or
 $$
 	    \beta -\alpha   \geq 2 \pi,
$$
with the solution confined ot the interval   $ [-1- \pi , 1] $ in the first case,
$ [-1, 2 \pi + 1 ] $ in the second case and   $ [\alpha , \beta ] $
in the last case. In particular, in the last case $ \beta - \alpha $
can be arbitrarily large.
}
\end{lem}

\vskip 0.2 in

\bcor \label{cor-pendulum} Consider the uncoupled system\footnote{Since we chose to concentrate on $ n = 4 $ pendula, we formulate the lemma for this case, although the proof 
carries over verbatim for an arbitrary $n$. } 
\begin{equation} 
	\ddot x_i+ \sin x_i=0, \ \ i=1, \ldots , 4. 
	\label{eq:uncoupled}\end{equation}   
For any $q, p \in {\mathbb R}  ^4 $  
such that $\al=q_i, \bt=p_i$ satisfies
conditions of Lemma \ref{lem:pendulumlemma} for each $i=1, \ldots , 4 $, then for any $T>0$
there is a unique solution $X (t,T,q,p)$ of  the system   (\ref{eq:uncoupled})   satisfying 
$X(0)=q$ and $X(T)=p$. Moreover,
there is a unique $T>0$ such that energy of this solution is one.
\ecor

We first prove the Corollary.

\vskip 0.2 in

\noindent{\bf Proof of the Corollary.}
 Lemma \ref{lem:pendulumlemma} applies to each of the   $n$ equations in   (\ref{eq:uncoupled}), by the assumptions of the Corollary.     
 That is, for any  $ T> 0 $
there exists a unique solution $ x_i(t) = x_i(t; T, q_i,p_i) $
of $ \ddot x_i + \sin x_i = 0$ with the desired boundary conditions.
Now as $T$ increases from $0$ to $ \infty $, the energy
$ E_i (T) $ of each solution decreases monotonically from $ \infty $
to $ 0 $. There is therefore a unique $T=T^\ast $ with
$   \Sigma_{i=1}^nE_i(T^\ast) = 1 $. In addition, $ T^\ast $ is
a smooth function of $ q , p $, by an application of the implicit
function theorem.   Q.E.D.

\vskip 0.3 in

\noindent{\bf Proof of the Lemma.} We concentrate on the first case; the remaining
ones are essentially the same. We write the equation of the pendulum as the  system
	\begin{equation}
   \left\{ \begin{array}{l}
    \dot x = y \\[3pt]
    \dot y  = -   \sin x . \end{array} \right.
 	\label{eq:pensys}
\end{equation}

\begin{figure}[thb]
         \psfrag{al}[c]{\small $  \alpha $}
 	 \psfrag{be}[c]{\small $ \beta $}
	  \psfrag{x}[c]{\small $x(t; T,\alpha , \beta)$}
	  \psfrag{p}[c]{\small $\pi $}
	  \psfrag{t}[c]{\small $t=0$}
	  \psfrag{T}[c]{\small $t=T$}
	  \psfrag{O}[c]{\small $O$}
	  \psfrag{M}[c]{\small $M$}
	  \psfrag{A}[c]{\small $A$}
	  \psfrag{B}[c]{\small $B$}
	  \psfrag{S}[c]{\small $S$}
	  \psfrag{IT}[c]{\small $I_T$}
	  \psfrag{I1}[c]{\small $I_{T_1}$}
	  \psfrag{T1}[c]{\small $ T_1\gg  T$}
    	\center{ \includegraphics[scale=0.7]{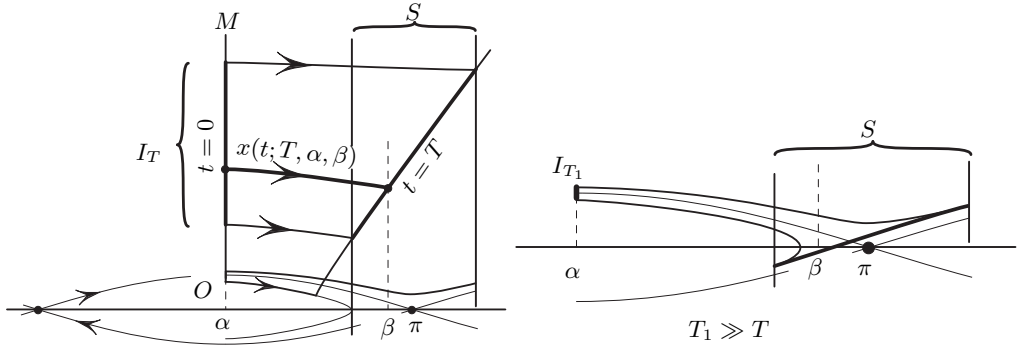}}
    	\caption{ The slope of
	  the image of the vertical interval is positive.}
    	\label{fig:pendulumlemma}
\end{figure}

We are seeking a solution starting on the line $ x= \alpha $ and ending at $ t=T $
on the line $ x= \beta $. The condition $ \beta \geq \pi -1 $ imposes a lower
bound on the initial velocity. This leads us to consider the ray $ OM $ of initial
data on the line $ x= \alpha $, see Figure ~\ref{fig:pendulumlemma}, where $O$ is
the point whose solution crosses the $x$--axis at $ x= \pi - 1 $.

 Consider now the ray $ OM $ carried by the flow of   (\ref{eq:pensys})
 for time $T$, where $T$ is fixed arbitrarily. Consider the set $ I_T $
 on the ray $ OM $  whose points enter  the strip
 $S=  \{ \pi -1 \leq x\leq \pi + 1\} $ at some time $ t  \leq T $ and
 do not leave it before $ t=T $.
   We claim: {\it  For any $T>0$, the set $ I_T $ is an interval, and
   its image $ \varphi ^T(I_T) $ under the flow at $ t=T $ is a curve
   with a positive slope, connecting the two boundaries of the strip $S$,
   Figure~\ref{fig:pendulumlemma}. } This claim implies the statement of the Lemma.

If $T$ is small, i.e. if the solutions are fast and thus lie above the separatrix,
the result is obvious. However,  for $T$ large  some   solutions
starting on $ I_T $ ``turn around", as in Figure~\ref{fig:pendulumlemma} (B),
and the proof requires a little care.
 \vskip 0.3 in
  Let $ z_0=( \alpha, y_0 ) \in I_T $; this, we recall, means that  (i)
  $z^T\equiv \varphi ^T z_0\in S $, and (ii) once  $ z^t $ enters $S$
  it does not leave $S$ before $ t=T $. Consider  the
linearization of (\ref{eq:pensys}):
\begin{equation}
   \left\{ \begin{array}{l}
   \dot \xi = \eta \\[3pt]
   \dot\eta= -(   \cos x )\; \xi ,  \end{array} \right.   \label{eq:linearized}
\end{equation}
where $ x$ is the solution of (\ref{eq:pensys}) with the chosen initial condition.
 Note that  the solution     $ \zeta = ( \xi, \eta) $ of (\ref{eq:linearized}) with
 $ \zeta(0)=(0, \eta_0),\ \eta_0>0$ is a tangent vector to the image curve
 $ \varphi ^T( OM) $ at the point $ \varphi ^Tz_0 $. It suffices, therefore,
 to prove that {\it $ \zeta(T) $ lies in the first
 quadrant. } This equation is changing type from elliptic to hyperbolic, and
 thus is it not {\it  a priori } clear that the solution may not execute an unwanted
 rotation, violating positivity of the slope $ \eta/\xi  $. The idea is to show
 that during the ``dangerous" elliptic stage while $ x(t)< \pi -1 $, $ \eta/\xi  $
 remains positive.

 To that end, let $ \tau \in (0,T] $ be the time of entrance of the solution
 $ \varphi ^t z_0=(x(t), y(t))$  into $S$, so that
 $$
 	x(t)\in [ \alpha , \pi -1 ]\ \ \hbox{for} \ \  t\in [0, \tau ]
$$
  and
  $$
 	 x(t) \in [\pi -1, \pi + 1]\ \ \hbox{for} \ \  t\in [  \tau, T ]
$$
 \vskip 0.2 in
  We will first show  that $ \eta/\xi >0 $ at $ t= \tau $.

The key idea is to   to observe that  the vector $ z= (x, \dot x ) $ rotates clockwise faster than the vector $\zeta=(\xi,\eta) $; this will be shown shortly. Since the slope of $z$ is positive at $ t= \tau $, the same will then be true of $ \zeta (\tau)$ as desired.
We claim:
\begin{equation}
	\frac{d}{dt} \biggl(\frac{y}{x}  \biggl) < \frac{d}{dt} \biggl(\frac{\eta}{\xi}  \biggl)< 0,
	\ \ \hbox{whenever}\ \ \frac{y}{x} =    \frac{\eta}{\xi}.
	\label{eq:key}
\end{equation}
To prove (\ref{eq:key}) we carry out the differentiations and use the equality of the slopes to reduce the inequality to an equivalent one:
  $$
  	\frac{\dot y  x - y \dot x }{ x ^2 } < \frac{\dot \eta  \xi - \eta \dot \xi }{ \xi ^2 } < 0 ,
  $$
  or, using (\ref{eq:pensys}) and (\ref{eq:linearized}):
 $$
 	\frac{-   x\sin   x - y^2  }{ x ^2 } <
	\frac{ -   \cos x \;   \xi ^2  - \eta ^2 }{ \xi ^2 } < 0.
 $$
  Using the equality of slopes, this reduces to
  $$
   	 \frac{\sin x}{x} > \cos x ,
  $$
   which holds true thanks to $ x\in (0, \pi ) $.
   We conclude: since $ y(\tau)/x(\tau) \geq 0 $, and $ \eta(\tau) / \xi (\tau) > 0 $.

 \vskip 0.2 in
  We show that the slope of $ \zeta $ remains positive for the remaining time
   $ [\tau,T]$. During this time $ | x- \pi | \leq 1 $ and thus $ \cos x < 0 $. Hence the linearized vector field
   (\ref{eq:linearized}) crosses {\bf into } the first quadrant, and since $ \zeta(\tau) $
   lies in that quadrant, it is still there at $ t=T$.
\vskip 0.2 in
 The monotonicity of $E=E(T)$ also follows from (i) the positivity of the slope of the
   image curve $ \varphi ^T(I_T) $, and (ii) from the fact that the curve moves ``to the right"
   as $T$ increases. Indeed, if we increase $T$, the point of intersection of $ \varphi ^T(I_T) $
   with $ x = \beta $ moves down, to the curve with the same $x= \beta $ but with
   the smaller $y$,
   i.e. with the smaller value of $ E= y ^2 +   (- 1- \cos \beta ) $. Furthermore,
   if $T$ is small,
   the velocity must be large: $\beta - \alpha=  x(T)-x(0) \leq ({\rm max} \; y)\; T $,
   so that $ {\rm max} \; y\geq \frac{\beta - \alpha}{T} \rightarrow \infty $ as
   $ T \rightarrow 0 $. On the other hand, if $ T $ is sufficiently large, then all
   solutions starting on $ I_T $ must pass arbitrarily close to the saddle $ ( \pi , 0) $
   and thus must have energy close to that of the saddle, i.e.  to
   $ E_{\rm saddle}=  - (- 1 - \cos \pi ) = 0   $.

\vskip 0.2 in
   The proof of Lemma \ref{lem:pendulumlemma} is complete.

 \section{The Hyperbolic Lemma.}\label{sec:hyperbolic}
 
By   Corollary~\ref{cor-pendulum} of Lemma~\ref{lem:pendulumlemma},  given any  points 
$ q, \ q ^\prime \in {\mathbb R} ^4  $ with their coordinates $ x_i $ and $ x ^\prime _i $ lying
 $ [-1,1]\, \hbox{mod} \; 2 \pi $ or in $  [ \pi -1,\pi +1] \ \hbox{mod} \; 2 \pi$, there exists a solution
 $ X(t; q, q ^\prime, T )$ of   (\ref{eq:uncoupled}) with $ \varepsilon = 0 $  which travels from $q$ to $ q ^\prime  $ in time $T$. By the same corollary,
 there exists a unique $T(q, q ^\prime ) $ for which the total energy of the solution
 $ X(t; q, q ^\prime, T(q, q ^\prime ))  $ is one:
 \be \label{eq:energy-one} 
 \sum_{i=1}^4 E_i =1, \ \ \hbox{where} \ \  E_i= \frac{\dot X_i ^2 }{2} + V(X_i).
 \ee
 We thus associate with the energy one solution  (of   (\ref{eq:mainsystem}) with $ \varepsilon = 0 $) connecting $ q $ and $ q ^\prime $,      the energy vector
 $$
 {\bf E} (q, q ^\prime ) \buildrel{def}\over{=} (E_1,E_2,E_3,E_4);
 $$
 according to the Lemma~\ref{lem:pendulumlemma}, this vector is uniquely determined by the endpoints 
 $ q , q ^\prime $. 
  \vskip 0.3 in
\begin{lem}\label{lem:velocity}   
If two pairs of points: $ q_1, q_1 ^\prime $ and $ q_2, q_2 ^\prime $ in $ {\mathbb R}  ^4 $ satisfy the conditions\footnote{Since we chose to concentrate 
on $ n = 4 $ pendula, we formulate the lemma for 
this case, although the proof 
carries over verbatim for an arbitrary $n$. } 
\begin{equation}
	| q_k ^\prime - q_k | \geq \varepsilon ^{-2r-5}, \ \ k= 1, 2,
	\label{eq:distancelarge}
\end{equation}
and
\begin{equation}
	| e_2-e_1 | < \varepsilon ^{2r+5},  \ \ \hbox{where} \ \
	e_k= \frac{ q_k ^\prime - q_k }{| q_k ^\prime - q_k |}, \ \ k =1, 2, 
	\label{eq:anglesmall}
\end{equation}
then the energy vector ${\bf E}$ of the connecting solution of      (\ref{eq:mainsystem}) with 
$ \beta = 0 $ satisfies
\begin{equation}
	| {\bf E} (q_2, q_2 ^\prime ) - {\bf E} (q_1, q_1 ^\prime )| < \varepsilon ^{2r+5} .
	\label{eq:energychange}
\end{equation}
	Moreover, there exists a constant $C$  such that for all $q$, $q ^\prime $ 
	with $ | q ^\prime - q | \geq 1 $ we have
\begin{equation}
	\biggl| \frac{d}{dq} \dot X(0; q, q ^\prime , T(q,q ^\prime )) \biggl| < C;
	\label{eq:boundedderivative}
\end{equation}
here the notation $ \frac{d}{dq}  $ is used to emphasize that the $q$--dependence enters $ \dot X $  in two places -- one through the boundary condition, and the other through
$T(q, q ^\prime ) $.
 \end{lem} 	

\noindent{\bf Proof.} Statement  (\ref{eq:energychange})  follows from the proof of 
Lemma~\ref{lem:pendulumlemma}; the main difficulty is in proving   (\ref{eq:boundedderivative})\footnote{this estimate can be strenghened: $ C $ can be replaced by $ C/|q ^\prime - q | $, but we do not need this in our proof.}.
To prove    (\ref{eq:boundedderivative}),we expand its left hand side:
\begin{equation}
 	\frac{d}{dq} \dot X(0; q, q ^\prime , T(q,q ^\prime ))=
	\partial_q  \dot X(0; q, q ^\prime , T ) +
	\partial_T \dot X(0; q, q ^\prime , T )\cdot \partial_q T(q, q ^\prime ),
	\label{eq:expanded}
\end{equation}
where $  T=T(q, q ^\prime )  $ is to be substituted after the differentiations on the right--hand side. We will now estimate each of the   summands on the right--hand side separately.
\paragraph {Estimate of $ \partial_q  \dot X(0; q, q ^\prime , T )$.} The proof of Lemma~\ref{lem:pendulumlemma} shows that each component $X_i$ of $X$ depends on the boundary conditions $ x_i, x_i ^\prime $ and $T$ only\footnote{We recall the notation $ q=(x_1,x_2,x_3,x_4) $.}, but not on     $ x_j, x_j ^\prime $ with $ j \not= i $. This implies that the matrix $ \partial_q  \dot X(0; q, q ^\prime , T )$ is diagonal, with the diagonal entries $  \partial X_i(x_i, x_i ^\prime , T)/ \partial x_i $.   But this derivative is simply the slope of the image   of the line $ x=x^\prime $ under the map $ \varphi ^{-T} $, where $\varphi^t $ is the phase flow of the pendulum equation.   The argument of Lemma~\ref{lem:pendulumlemma} shows that, because of the shear in the phase velocity field, this slope is always bounded once $T$ exceeds a fixed constant. It remains to prove the upper bound for the last summand in   (\ref{eq:expanded}).

\paragraph {Estimate of $ \partial_T \dot X(0; q, q ^\prime , T )\cdot \partial_q T(q, q ^\prime ) $.}
We will first show that this term is expressible via the first factor alone, thus reducing the number of estimates needed. Note that   $ \partial_T \dot X \cdot \partial_q T $ is a square matrix, the product  of the row matrix $ \partial_T \dot X $ with the column gradient matrix $ \partial_q T(q, q ^\prime ) $. To prove the lemma, it remains to show that each entry
\begin{equation}
	| \partial_T \dot X_i(0; x_i, x_i  ^\prime , T )\cdot \partial_{x_j} T(q, q ^\prime ) | < C,
	\label{eq:boundedentry}
\end{equation}
for $ | q ^\prime - q | \geq 1 $.

\paragraph {Some identities.} Let
\begin{equation}
	K_i= \biggl( \int_{x_i}^{x_i ^\prime } \frac{dx}{(2(E_i-V(x))) ^{3/2}  } \biggl)^{-1}, \ \ \hbox{and} \ \ K = 	\sum_{s=1}^4 K_s.
	\label{eq:K}
\end{equation}
We will show that
\begin{equation}
	\partial_{x_i} T(q, q ^\prime ) = K ^{-1}  \partial_T \dot X_i(0; x_i, x_i  ^\prime , T ),
	\label{eq:relation}
\end{equation}
thus reducing     (\ref{eq:boundedentry})  to an equivalent inequality
\begin{equation}
	| K ^{-1}\;  \partial_T \dot X_i(0; x_i, x_i  ^\prime , T )\; \partial_T \dot X_j(0; x_j, x_j  ^\prime , T ) | < C.
	\label{eq:boundedentry1}
\end{equation}
Heuristically, one expects that   $ | \partial_T \dot X_i(0; x_i, x_i  ^\prime , T ) | \leq c T ^{-1}  $.
Indeed,    $   \dot X_i(0; x_i, x_i  ^\prime , T ) $ is the $y$--coordinate of the intersection, in the
$ (X, \dot X)$--plane  of the line $ \{X=x_i\} $ and the curve $ \varphi ^{-T}  \{X=x_i^\prime \} $,
where $ \varphi ^t $ is the phase flow of the pendulum equation. Now because of the shear in the phase flow, one expects the line $ \ell_T = \varphi ^{-T}  \{X=x_i^\prime \} $ to form angle at most
$ c T ^{-1} $ with the trajectories. Thus the point $ \ell_T\cap \{X=x_i\} $ is expected to move with speed
$ \leq c T ^{-1} $, suggesting that indeed $ | \partial_T \dot X_i(0; x_i, x_i  ^\prime , T ) | \leq c T ^{-1}  $.
We carry out a precise proof by an alternative, purely analytical method (which ultimately reduces to the same estimates). Namely, we will use the following identity:

\begin{equation}
	\partial_T \dot X_i(0; x_i, x_i  ^\prime , T ) = -\frac{K_i}{ \sqrt{2( E_i-V(x_i))}},
	\label{eq:xdotbyT}
\end{equation}
which, together with   (\ref{eq:relation}), is  proven in a separate section below.
\vskip 0.3 in
 
\paragraph {Estimate of $ K ^{-1} $.} Since $ \Sigma_{i=1}^4 E_i = 1 $, we have $ \frac{1}{4}\leq E_i\leq 1 $ for some $i$, and thus for some $C$ we have
$$
	K_i ^{-1} =  \int_{x_i}^{x_i ^\prime } \frac{dx}{ (2( E_i-V(x) )) ^{3/2} }  \leq	
	C\int_{x_i}^{x_i ^\prime } \frac{dx}{  \sqrt{ 2( E_i-V(x)) }  }    = C  T ,
$$
so that
\begin{equation}
	K ^{-1} = \biggl( \sum_{j=1}^4 K_j \biggl) ^{-1} <  K_i ^{-1} \leq CT.
	\label{eq:Kestimate}
\end{equation}

\paragraph {Estimate of  $ \partial_T \dot X_i(0; x_i, x_i  ^\prime , T )$.} We     consider two separate cases: (i) $ |x_i ^\prime -x_i | < 2 \pi $ and (ii) $ |x_i ^\prime -x_i | \geq 2 \pi $.

\noindent 1. {\bf In  case (i)}, an estimate of    (\ref{eq:xdotbyT}) is easier done geometrically, as follows. Consider
the graph $ y = U_T(x) $ of the time $T$--preimage of the line $x= x'$
in the phase plane of the pendulum. Let $ y = U(x)   $ be the graph of the unstable manifold of the saddle
$ ( \pi , 0) $;   by a standard hyperbolic argument, the flow takes the line exponentially close to the unstable manifold:
 $ | U_T(x)-U(x) | <e^{-cT}  , \ \ {\rm for} \ \ |x| \leq \pi $,
and, moreover, the motion of the line becomes exponentially slow:
\begin{equation}
	|\frac{d}{dT}  U_T(x)| <e^{-cT}  , \ \ {\rm for} \ \ |x| \leq \pi;
	\label{eq:xdotbyT1}
\end{equation}
 here $T$ is greater than a fixed positive constant because of the assumption
$ | q ^\prime - q | \geq 1 $. But $ U_T(x_i)= \dot X_i(0; q, q ^\prime, T)$ and     (\ref{eq:xdotbyT1}) gives
\begin{equation}
	| \partial_T \dot X_i(0; q, q ^\prime, T)| \leq e^{-cT}, \ \ \hbox{for} \ \ |x_i ^\prime -x_i | < 2 \pi.
	\label{eq:xdotbyTexp}
\end{equation}
This completes the proof of   (\ref{eq:boundedentry1}),  and thus of the Lemma, in case (i).

\noindent 2. {\bf In case (ii)}  we have $   x_i ^\prime - x_i   = 2 \pi n_i + r, \ \ 0\leq r< 2 \pi $ with integer $n\not=0$. We will use   (\ref{eq:xdotbyT}) to prove   (\ref{eq:boundedentry1}), to which end we       need an upper bound on $ K_i $. From   (\ref{eq:K})  we have
\begin{equation}
	K_i^{-1} \geq   n_i \int_{- \pi }^{  \pi } \frac{dx}{ (E_i-V(x))^{3/2}}  \geq
	n_i \int_{- \pi }^{  \pi } \frac{dx}{ (E_i+x ^2 /2)^{3/2}} = c \frac{n_i}{E_i},
	\label{eq:Kibound}
\end{equation}
where $ c= 2 \pi / \sqrt{ 1+ \pi ^2 /2 }    $.
Now the number of revolutions $ n _i \geq T/{\cal T} _{E_i}$ where  ${\cal T}(E_i) $  is  the time of one full revolution:
as
$$
	{\cal T} _{E_i}= \int_{- \pi }^{  \pi } \frac{dx}{ \sqrt{ 2(E_i+2 \sin ^2 x/2)} }\leq
	 \sqrt{ 2} \int_{ 0 }^{  \pi } \frac{dx}{ \sqrt{  E_i+x ^2 / \pi ^2 } } =
	 \sqrt{ 2 } \pi (\ln (1+  \sqrt{ 1+E_1})-\ln E_i).
$$
  Substituting this into   (\ref{eq:Kibound}) we get
$$
  	K_i\leq  c\frac{E_i}{n_i} \leq c\frac{E_i {\cal T_{E_i}} }{T} \leq
	c_1\frac{E_i (\ln (1+  \sqrt{ 1+E_i})-\ln E_i) }{T}  .
$$
Finally, we substitute this estimate into   (\ref{eq:xdotbyT}):
$$
	| \partial_T \dot X_i(0; x_i, x_i  ^\prime , T )| \leq
	c \frac{E_i (\ln (1+  \sqrt{ 1+E_i})-\ln E_i) }{T \sqrt{ E_i } } \leq \frac{c}{T}
$$
Together with   (\ref{eq:Kestimate}) this  proves   (\ref{eq:boundedentry1}). The proof of the lemma is thus complete.
\vskip 0.3 in

\paragraph {Proof of the identities   (\ref{eq:relation})  and   (\ref{eq:xdotbyT}).}
The energy of the solution $ X(t; x_i,x_i ^\prime , T) $ is a smooth function of Let $  x_i,x_i ^\prime , T  $; we denote this energy by   $ E(x_i,x_i ^\prime , T) $, so that
$$
\dot X(0; x_i, x_i ^\prime , T) =  \sqrt{2( E(x_i,x_i ^\prime , T)-V(x_i))}.
$$
Differentiating by $T$ we get
\begin{equation}
	\frac{ \partial}{\partial T}  \dot X(0; x_i, x_i ^\prime , T) =
	\frac{ \partial E(x_i,x_i ^\prime , T)/ \partial T}{\sqrt{2( E(x_i,x_i ^\prime , T)-V(x_i))}}
	\label{eq:ddoxdt}
\end{equation}
To estimate the numerator, we differentiate the identity
\begin{equation}
	T= \int_{x_i}^{x_i ^\prime } \frac{dx}{\sqrt{2( E(x_i,x_i ^\prime , T)-V(x))}}
		\label{eq:T}
\end{equation}
with respect to $T$ and solve for $ \partial E/\partial T $, obtaining
\begin{equation}
	\partial E(x_i,x_i ^\prime , T)/ \partial T = - \biggl(  \int_{x_i}^{x_i ^\prime }
	 \frac{dx}{ (2( E(x_i,x_i ^\prime , T)-V(x))) ^{3/2} }\biggl) ^{-1} \buildrel{def}\over{\equiv} - K_i ^{-1} .
 	\label{eq:ebyt}
\end{equation}
Substituting this into     (\ref{eq:ddoxdt})  proves   (\ref{eq:xdotbyT}).
To prove the remaining identity (\ref{eq:relation}), we recall that  $T(x, x ^\prime ) $ is the time which gives energy one to the solution:
\begin{equation}
	\sum_{k=1}^4 E(x_k, x_k ^\prime , T(x_k, x_k ^\prime)) = 1.
	\label{eq:energyone}
\end{equation}
Differentiating this by $ x_i $ gives
\begin{equation}
 	\frac{\partial E(x_i, x_i ^\prime , T)}{\partial x_i} \biggl|_{T=T(x , x  ^\prime)}\ + \
	 \frac{\partial T(x , x  ^\prime)}{\partial x_i} \ \sum_{k=1}^4\frac{\partial E(x_k, x_k ^\prime , T)}	
	 {\partial T}=0.
		 \label{eq:energyone1}
\end{equation}
The above sum, according to   (\ref{eq:ebyt}), can be replaced by
 $- \sum_{k=1}^4K_k \buildrel{def}\over{=}  -K$; solving      for $ \partial T/\partial x_i $ gives
 \begin{equation}
	\frac{\partial T(x , x  ^\prime)}{\partial x_i} =
	K ^{-1}\frac{\partial E(x_i, x_i ^\prime , T)}{\partial x_i} \biggl|_{T=T(x , x  ^\prime)}.
		\label{eq:tbyx}
\end{equation}
To estimate the last derivative, we differentiate the
identity\footnote{Here and below $ E_i=E(x_i,x_i ^\prime , T) $.}   (\ref{eq:T})  by $ x_i $:
$$
	0= -\frac{1}{\sqrt{2( E_i-V(x_i))}}-
	\underbrace{\int_{x_i}^{x_i ^\prime } \frac{dx}{ (2( E_i-V(x))) ^{3/2} }}_{K_i ^{-1} }\;
	\frac{\partial E(x_i, x_i ^\prime , T)}{\partial x_i},
$$
 or
$$
	\frac{\partial E(x_i, x_i ^\prime , T)}{\partial x_i} = - \frac{K_i}{\sqrt{2( E_i-V(x_i))}} .
$$
Substituting   this into   (\ref{eq:tbyx})  results in the proof of   (\ref{eq:relation}):
$$
 \frac{\partial T(x , x  ^\prime)}{\partial x_i} = -K ^{-1} \frac{K_i}{\sqrt{2( E_i-V(x_i))}}
 \buildrel{ (\ref{eq:xdotbyT})}\over{=} K ^{-1} \partial_T \dot X_i(0; x_i, x_i  ^\prime , T ).
$$
The proof of the two identities is now complete.


\section{The connection Lemma. } \label{sec:existence}
The following lemma is the building block in the construction of shadowing geodesics.

\begin{lem}[Existence of geodesic segments ]\label{lem:existence}\footnote{Since we chose 
to concentrate on $ n = 4 $ pendula, we formulate 
the lemma for this case, although the proof 
carries over verbatim for an arbitrary $n$. } There exists $ \varepsilon_0 >0 $ such that for all $0< \varepsilon < \varepsilon_0 $ the following holds.  Consider any two sections, which we denote by  $\Sigma_0$ and $\Sigma_1$, in the itinerary  (\ref{eq:sequencofballs}), such that the vector $2\pi \vec n$ connecting the centers of  $\Sigma_0, \ \Sigma_1$ satisfies $|\vec n|>\frac 1\eps$. This vector is  of the form 
$2\pi(m,n,1,0)$  or $2\pi (\frac 12 , s, 1, \frac 12)$ (with
  integer $m,n$ or $s$), modulo possible translation in the index.
Then for all $ p_0\in \Sigma_0 $, $ p_1\in \Sigma_1 $ there exists a geodesic 
$ \gamma (p_0, p_1) $ in the Jacobi metric    (\ref{eq:jacobimetric})      connecting 
$p_0$ with $ p_1$, and depending smoothly on $ p_0 $ and $  p_1  $.
\end{lem}


 \paragraph {Proof.} 
      \begin{figure}[thb]
         \psfrag{p0}[c]{\small $p_0$}
         \psfrag{p1}[c]{\small $p_1$}
 	 \psfrag{q0}[c]{\small $q_0$}
	 \psfrag{q1}[c]{\small $q_1$}
	  \psfrag{X0}[c]{\small $X_0(t)$}
	  \psfrag{c0}[c]{\small $c_0$}
	  \psfrag{c1}[c]{\small $c_1$}
	  \psfrag{t0}[c]{\small $t_0$}
	  \psfrag{t1}[c]{\small $t_1$}
	  \psfrag{vl}[c]{\small $v_L$}
	  \psfrag{vr}[c]{\small $v_R$}
	  \psfrag{wl}[c]{\small $w_L$}
	   \psfrag{wr}[c]{\small $w_R$}
	   \psfrag{S0}[c]{\small $S_0$}
	   \psfrag{S1}[c]{\small $S_1$}
    	\center{ \includegraphics[scale=0.7]{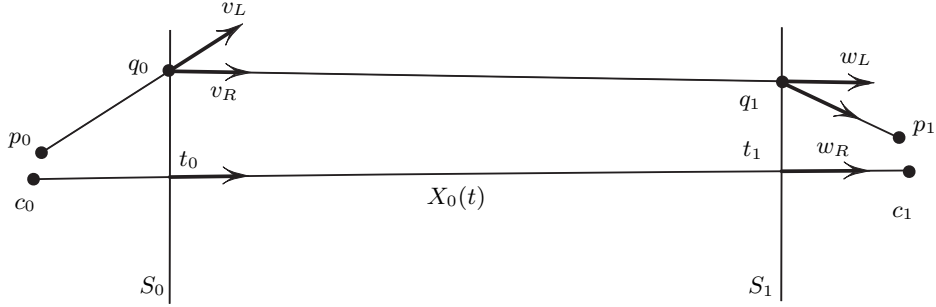}}
    	\caption{Towards proof of Lemma \ref{lem:existence}. }
    	\label{fig:transversal}
\end{figure}
 \begin{enumerate}

\item  We first define a section $  S_0$  shown in 
 Figure~\ref{fig:transversal}, as follows. By Lemma~\ref{cor-pendulum} there exists a  (unique) solution $ X_0(t)$ of   (\ref{eq:uncoupled}) of energy one,
 connecting $c_0=Center(\Sigma_0)$ with $c_1=Center(\Sigma_1)$.  
 We now define the section $\mathcal S_0$ as the codimension one disk in 
 ${\mathbb R}  ^4 $
 at the distance $ \varepsilon ^{\frac{1}{3} } $ from $ c_0 $ and perpendicular to the 
 initial direction $ e_0= \dot X_0(0)/|\dot X_0(0)| $:  
 $$
 	S_0= \{ q: (q-c_0)\cdot e_0 = \varepsilon^ \frac{1}{3}, \ \ 
 	| q-(c_0+ \varepsilon ^ \frac{1}{3} e_0) | < \varepsilon ^ \frac{1}{3} \},
 $$
 where $ \cdot $ denotes the usual dot product. 
Analogously, we define the section $ S_1 $ near $ c_1 $, except for reversing the sign of 
$ e_1$.  
 
 \item 
 For any pair of points $ q_0\in S_0 $, $ q_1\in S_1 $ we consider three geodesic 
 segments (in the metric   (\ref{eq:jacobimetric})):    $ \gamma (p_0, q_0) $, 
 $  \gamma (q_0, q_1) $ and  $ \gamma (q_1, p_1) $, along with    the velocities of the associated solutions of   (\ref{eq:mainsystem})
 $ v_L, \ v_R, \ w_L, w_R $ as shown in 
 Figure~\ref{fig:transversal}\footnote{It should be noted that $ v_L, \ v_R, \ w_L, w_R $ all depend on $ q_0, \ q_1 $.}.  The lemma will be proven
 once we show that there exists a pair $ q_0, \ q_1 $, smoothly dependent on $ p_0, \ p_1 $, 
  for which 
 \begin{equation} 
	 v_L= v_R \ \ \hbox{and }\ \  w_L = w_R.
	\label{eq:matching}\end{equation}   
To that end we first list the properties of each of the three geodesic segments.  

\item   Since radius of $\Sigma_0$ is $\eps^{1/2}$, we have  
$ | q_0-p_0 | = O (\varepsilon ^{ \frac{1}{3} } )$, which  is small compared to 
the injectivity radius of the metric   (\ref{eq:jacobimetric}).   
By standard arguments from differential geometry (using the smoothness 
of solutions of the ODEs and the implicit function theorem) we conclude that 
 \begin{equation} 
	 v_L= v_0 \frac{q_0-p_0}{|q_0-p_0|} + r_{0L}(p_0,q_0, \varepsilon ),  \ \  
	 | r_{0L}|_{C^1} = O( \varepsilon^\frac{1}{3} );  
	\label{eq:vl}\end{equation}   
here $ v_0= 3  \sqrt{ 2 } $    is the speed at $ c_0$ (by the energy 
restriction (\ref{eq:energy-one} and the fact that the center of the lens corresponds to three pendula at the bottom, and one at the top.). 
A similar estimate holds for the right end: 
\begin{equation} 
	 w_R= v_0\frac{p_1-q_1}{|p_1-q_1|} + r_{1R}(q_1,p_1, \varepsilon ),  \ \  
	 | r_{1R}|_{C^1} = O( \varepsilon^\frac{1}{3} ) . 
	\label{eq:wr}\end{equation}   

\item  The intermediate segment $ \gamma (q_0, q_1 ) $ avoids the lenses, and thus     Lemma 
\ref{lem:velocity} applies; in particular, the $ C^1 $--bound   (\ref{eq:boundedderivative}) holds,  
implying that 
\begin{equation} 
	v_R= \dot X_0(t_0)+ r_{0R}(q_0,q_1, \varepsilon ), \ \ 
	| r_{0R} |_{C^1}< C \varepsilon ^ \frac{1}{3} . 
	\label{eq:vrwl}\end{equation}   
and 
\begin{equation} 
	w_L= \dot X_0(t_1)+ r_{1L}(q_1,p_1, \varepsilon ), \ \ 
	| r_{1L} |_{C^1}< C \varepsilon ^ \frac{1}{3} . 
	\label{eq:vrwl1}\end{equation}
Here $ t_i $ ($i=0,1$) is  the time  when $ X_0 (t) $ intersects the section $ S_i$. 
\item We will prove the existence of the pair $ q_0, \ q_1 $ satisfying   (\ref{eq:matching})
by applying the implicit function theorem. To that end let $ \widehat v $ denote the 
orthogonal projection of 
$ v\in {\mathbb R}  ^4 $ onto ${\mathbb R}  ^3 \supset  S$; we will also treat $ q_0 \in {\mathbb R}  ^4 $ as an
 element of $ S\subset {\mathbb R}  ^3 $, denoting it  by $ \widehat q_0 \in {\mathbb R}  ^3 $.  
To prove   (\ref{eq:matching})  it suffices to prove that the projected equations
 \begin{equation} 
	 \widehat{v}_L = \widehat{v}_R \ \ \hbox{and }\ \  \widehat{w}_L = \widehat{w}_R 
	\label{eq:matching1}
\end{equation}  
hold, since if   (\ref{eq:matching1})  hold, then the remaining components orthogonal to $S$ must match as well by the conservation of energy. 
Substituting the  estimates   (\ref{eq:vl}),   (\ref{eq:wr}),    (\ref{eq:vrwl})  and   
(\ref{eq:vrwl1}) into   (\ref{eq:matching})  we obtain, after projecting onto 
$ {\mathbb R}  ^3 \supset S $: 
\begin{equation} 
	  v_0 \frac{\widehat{q_0-p_0}}{|q_L-p_0|}  = 
	  \widehat{r}_0(p_0,\widehat{q}_0,\widehat{q}_1, \varepsilon ), \ \ 
	 v_0\frac{\widehat{p_1- q_1}}{|p_1-q_1|} = 
	 \widehat{r}_1(\widehat{q}_0,\widehat{q}_1,p_1, \varepsilon )
	\label{eq:matching2}\end{equation}   
as the equivalent matching conditions, with the $ C^1 $--small remainders: 
$$
| \widehat{r}_i|_{C^1}< C \varepsilon ^ \frac{1}{3}.  
$$   
In arriving at   (\ref{eq:matching2}), we made use of the fact that $ \widehat {\dot X_0(t_0)}= 
O( \varepsilon ^ \frac{1}{3} )$, 
as follows from    the choice of $ S_0$ to be orthogonal to $ \dot X_0(0) $ (so that 
$ \widehat{\dot X_0(0) } =0$) and the the fact that $ t_0 = O(\varepsilon ^ \frac{1}{3} )$.  
 \item To apply the implicit function theorem, instead of the variables $ \widehat{q}_i $ 
 we introduce
 $$
 	Q_0= v_0 \frac{\widehat{q_0-p_0}}{|q_L-p_0|} , \ \ 
	Q_1 = v_0\frac{\widehat{p_0-q_1}}{|p_1-q_1|},  
 $$
$ Q_i\in {\mathbb R}  ^3 $. Expressing 
$$
	\widehat{q}_0 = \widehat{p}_0+ v_0|q_0-p_0|Q_0, \ \ 
	\widehat{q}_1 = \widehat{p}_1- v_0|p_1-q_1| Q_1 
$$
 and substituting into 
(\ref{eq:matching1}), we obtain 
$$
Q_0 = R_0(p_0, Q_0, Q_1, \varepsilon ) , \ \ Q_1 = R_1(Q_0, Q_1, p_1\varepsilon ). 
$$
Introducing $ Q=(Q_0, Q_1) \in {\mathbb R}  ^6 $ and $ R=(R_0, R_1) $ we rewrite the matching condition   (\ref{eq:matching})  in the final form 
$$
Q= R(Q, p_0, p_1, \varepsilon ), 	
$$
where $ | R | _{C^1}< C \varepsilon ^ \frac{1}{3}$.  
It is important to observe that  $R$ is defined (at least) on the entire  ball 
$ | Q | \leq \frac{1}{2}   $, independent of $\varepsilon$. Thus for all sufficiently small 
$\varepsilon$ there exists a unique solution $ Q $ depending differentiably 
 on the parameters $ p_0, \; p_1 $. 
 \end{enumerate} 
This completes the proof of   Lemma \ref{lem:existence}.   

\end{document}